\documentclass[10pt]{article}
\usepackage{amssymb}
\usepackage{amsmath}
\usepackage{epsfig}

\normalsize

\newtheorem{theorem}{Theorem}[section]
\newtheorem{lemma}[theorem]{Lemma}

\newtheorem{corollary}[theorem]{Corollary}
\newtheorem{algorithm}[theorem]{Algorithm}

\newtheorem{example}[theorem]{Example}

\newcommand{\boproof}{\noindent {\bf Proof. }}
\newcommand{\eoproof}{\hspace*{\fill} $\square$ \vspace{5pt}}
\newcommand{\red}{\sqsubseteq}

\newcommand{\SA}{{{\cal S}_\Lambda}}
\newcommand{\R}{\mathbb R}

\newcommand{\Z}{\mathbb Z}

\newcommand{\GIP}{{\cal G}}
\newcommand{\symGIP}{{\cal G}^{\text{sym}}}

\newcommand{\Kw}[1]{\underline{#1}}

\DeclareMathOperator{\orb}{orb_{S_\Lambda}}
\DeclareMathOperator{\rep}{rep_{S_\Lambda}}
\DeclareMathOperator{\nF}{normalForm}
\DeclareMathOperator{\4ti2}{4ti2}

\begin{document}
\title{Exploiting Symmetries in the Computation of Graver Bases}
\author{Raymond Hemmecke\\Otto-von-Guericke-University Magdeburg\\raymond@hemmecke.de}
\date{}
\maketitle

\begin{abstract}
Many challenging Graver bases computations, like for multi-way
tables in statistics, have a highly symmetric problem structure
that is not exploited so far computationally. In this paper we
present a Graver basis algorithm for sublattices of $\Z^n$ that
exploits existing symmetry.
\end{abstract}

%\paragraph{Key Words.}
%\paragraph{AMS subject classifications.}

\section{Introduction} \label{Introduction}

Graver bases, originally introduced by Graver \cite{Graver:75} for
use in integer programming, have a variety of interesting
applications. Besides providing improving directions for integer
programs \cite{Graver:75,Hemmecke:PSP,Weismantel:98}, for
stochastic integer programs \cite{Hemmecke:SIP2}, and even for
certain convex integer programs
\cite{Hemmecke:Zconvex,Murota+Saito+Weismantel:04}, Graver bases
can for example also be used as Markov bases for sampling in
statistics \cite{Diaconis+Sturmfels:98}, or as a superset from
which a universal Gr{\"o}bner basis of the toric ideal
$I_A:=\langle x^u-x^v:Au=Av,u,v\in\Z^n_+\rangle$ can be extracted
\cite{Sturmfels:95}.

Unfortunately, the size of Graver bases increases quickly with the
dimension, making it very hard if not impossible to compute them
in practice. In several applications, however, as for example in
algebraic statistics, the problems involve a high symmetry that
should make it much easier to compute the Graver basis in terms of
(relatively few) representatives. In
\cite{Jach:DA,Jach+Koeppe+Weismantel:04}, the authors exploit
existing symmetry of lattices that arise from a single constraint
of the form $a^\intercal x\equiv 0\pmod p$, $p\in\Z_+$, in a
similar way as we do in Lemma \ref{Orbits in G(A) are always full}
and in Corollary \ref{Cor: Orbits in G(A) are always full} below.

Let us start our presentation by defining the notion of a Graver
basis and by giving an example that demonstrates the problem we
are interested in.

The {\bf Graver basis} $\GIP(\Lambda)$ associated to a lattice
$\Lambda\subseteq\Z^n$ consists exactly of all $\red$-minimal {\it
nonzero} elements in $\Lambda$, where for $u,v\in\Z^n$ we say that
$u\red v$ if $u^{(j)}v^{(j)}\geq 0$ and $|u^{(j)}|\leq |v^{(j)}|$
for all components $j=1,\ldots,n$, that is, if $u$ belongs to the
same orthant as $v$ and its components are not greater in absolute
value than the corresponding components of $v$. Note that Graver
originally defined this set only for the case
$\Lambda=\ker(A)\cap\Z^n$ for given matrix $A\in\Z^{d\times n}$,
but his definition can be readily extended to the definition we
gave above.

\begin{example} \label{Example: 3x3 tables} {\rm Consider the set of all
$3\times 3$ tables/arrays whose entries are filled with integer
numbers in such a way that the sums along each row and along each
column are $0$. One particular example is the table
\[
\left(
\begin{array}{rrr}
 1 & -1 &  0\\
-1 &  3 & -2\\
 0 & -2 &  2\\
\end{array}
\right).
\]
If we encode the $9$ entries of the table as $z_1,\ldots,z_9$,
then the set of $3\times 3$ tables coincides with the integer
vectors in the kernel of the matrix
\[
A_{3\times 3}=\left(
\begin{array}{rrrrrrrrr}
1 & 1 & 1 & 0 & 0 & 0 & 0 & 0 & 0\\
0 & 0 & 0 & 1 & 1 & 1 & 0 & 0 & 0\\
0 & 0 & 0 & 0 & 0 & 0 & 1 & 1 & 1\\
1 & 0 & 0 & 1 & 0 & 0 & 1 & 0 & 0\\
0 & 1 & 0 & 0 & 1 & 0 & 0 & 1 & 0\\
0 & 0 & 1 & 0 & 0 & 1 & 0 & 0 & 1\\
\end{array}
\right),
\]
that is, with all $z\in\Z^9$ satisfying $A_{3\times 3}z=0$. As
defined above, the Graver basis of $A_{3\times 3}$ consists of all
$\red$-minimal {\it nonzero} tables among them. The particular
$3\times 3$ table above does not belong to the Graver basis of
$A_{3\times 3}$, since
\[
\left(
\begin{array}{rrr}
 1 & -1 &  0\\
-1 &  1 &  0\\
 0 &  0 &  0\\
\end{array}
\right)\red
\left(
\begin{array}{rrr}
 1 & -1 &  0\\
-1 &  3 & -2\\
 0 & -2 &  2\\
\end{array}
\right).
\]

Using the computer program $\4ti2$ \cite{4ti2}, we find that the
following $15$ vectors (and their negatives) constitute the Graver
basis of $A_{3\times 3}$:
\[
\begin{array}{rrrrrrrrr}
(1, & -1, &  0, & -1, &  1, &  0, &  0, &  0, &  0)\\
(0, &  0, &  0, &  1, &  0, & -1, & -1, &  0, &  1)\\
(1, &  0, & -1, & -1, &  0, &  1, &  0, &  0, &  0)\\
(1, & -1, &  0, &  0, &  0, &  0, & -1, &  1, &  0)\\
(0, &  0, &  0, &  1, & -1, &  0, & -1, &  1, &  0)\\
(1, & -1, &  0, & -1, &  0, &  1, &  0, &  1, & -1)\\
(0, & -1, &  1, &  1, &  0, & -1, & -1, &  1, &  0)\\
(1, & -1, &  0, &  0, &  1, & -1, & -1, &  0, &  1)\\
(1, &  0, & -1, &  0, &  0, &  0, & -1, &  0, &  1)\\
(0, & -1, &  1, &  0, &  1, & -1, &  0, &  0, &  0)\\
(0, &  1, & -1, &  1, & -1, &  0, & -1, &  0, &  1)\\
(0, &  0, &  0, &  0, &  1, & -1, &  0, & -1, &  1)\\
(0, & -1, &  1, &  0, &  0, &  0, &  0, &  1, & -1)\\
(1, &  0, & -1, &  0, & -1, &  1, & -1, &  1, &  0)\\
(1, &  0, & -1, & -1, &  1, &  0, &  0, & -1, &  1)\\
\end{array}
\]
However, there is an obvious symmetry group $S_3\times S_3\times
S_2$ operating on the set of $3\times 3$ tables whose elements
transform a given table $v\in\ker(A_{3\times 3})$ into another
table $w\in\ker(A_{3\times 3})$ by suitably rearranging components
(permuting rows or columns, flipping the table along the main
diagonals). If we take these symmetries into account, we see that
among these 15 elements there are in fact only two essentially
different elements:
\[
\begin{array}{rrrrrrrrr}
(1, & -1, &  0, & -1, &  1, &  0, &  0, &  0, &  0)\\
(1, & -1, &  0, & -1, &  0, &  1, &  0, &  1, & -1)\\
\end{array}
\]
or, in a more array-like notation:
\[
\left(
\begin{array}{rrr}
 1 & -1 &  0\\
-1 &  1 &  0\\
 0 &  0 &  0\\
\end{array}
\right)
\;\;\;\text{and}\;\;\;
\left(
\begin{array}{rrr}
 1 & -1 &  0\\
-1 &  0 &  1\\
 0 &  1 & -1\\
\end{array}
\right).
\]
It should be clear that for bigger or for higher-dimensional
tables, this difference in sizes becomes far more striking, since
the acting symmetry groups are much bigger.} \eoproof
\end{example}

In the following, let a sublattice $\Lambda\subseteq\Z^n$ be given
and let $\SA\subseteq S_n$ be a group of symmetries such that for
all $v\in\Lambda$ and for all $\sigma\in\SA$ we have that also
$\sigma(v):=(v^{\sigma(1)},\ldots,v^{\sigma(n)})\in\Lambda$.
Finally, denote by $\orb(v):=\{\sigma(v):\sigma\in\SA\}$ the orbit
of $v$ under $\SA$.
\begin{lemma}\label{Orbits in G(A) are always full}
Let $\Lambda\subseteq\Z^n$ be a lattice, let $\SA\subseteq S_n$ be
the group of its symmetries, and let $g,g',s,s'\in\Lambda$. Then
the following holds:
\begin{itemize}
\item If $g\red s$ then $\sigma(g)\red\sigma(s)$ for every
$\sigma\in\SA$.
\item If there is no $g'\in\orb(g)$ with $g'\red s$, then for every
$s'\in\orb(s)$ there is no $g''\in\orb(g)$ with $g''\red s'$.
\item If $v\in\GIP(\Lambda)$ then $\orb(v)\subseteq\GIP(\Lambda)$.
\end{itemize}
\end{lemma}
\boproof The first statement follows immediately from the
definition of $\red$.

For the second statement, let $s'=\sigma(s)$ for some
$\sigma\in\SA$ and assume there is some $g''\in\orb(g)$ with
$g''\red s'$. Then $g':=\sigma^{-1}(g'')\red\sigma^{-1}(s')=s$ and
$g'\in\orb(g)$. A contradiction to the assumed non-existence of
such $g'$.

For the last statement we have to show that with
$v\in\GIP(\Lambda)$ also the full orbit $\orb(v)$ lies in
$\GIP(\Lambda)$. For this it suffices to assume that there is some
$\sigma(v)\in\orb(v)$ that could be written non-trivially as
$\sigma(v)=w_1+w_2$ with $w_1,w_2\red\sigma(v)$. But this would
imply $v=\sigma^{-1}(w_1)+\sigma^{-1}(w_2)$ with nonzero
$\sigma^{-1}(w_1),\sigma^{-1}(w_2)\red v$, which contradicts
$v\in\GIP(\Lambda)$. Thus, $\sigma(v)$ must belong to
$\GIP(\Lambda)$. \eoproof

As a consequence of this lemma, $\GIP(\Lambda)$ decomposes
completely into {\it full} orbits. Our task of computing
$\GIP(\Lambda)$ thus reduces to computing representatives of these
orbits, and to collect them into a set $\symGIP(\Lambda)$. By the
previous lemma, we recover $\GIP(\Lambda)$ via
\[
\GIP(\Lambda)=\bigcup_{v\in\symGIP(\Lambda)}\orb(v).
\]
Note that this last expression does not compute a superset of
$\GIP(\Lambda)$. In contrast to this, the last statement of Lemma
\ref{Orbits in G(A) are always full} fails to be true in general
for minimal toric Gr{\"o}bner bases or for minimal Markov bases
associated with the lattice $\Lambda$.

\section{Computing the Graver Basis}
\label{Section: Computing the Graver Basis}

In this section, we adapt Pottier's algorithm \cite{Pottier:96} to
deal with the symmetries of $\Lambda$ in the computation of
$\GIP(\Lambda)$. Note that this algorithm is not the fastest way
to compute Graver bases directly. However, by adapting this
algorithm, it will be easier for us to exploit the given
symmetries (and to present the main ideas). The state-of-the-art
algorithm that is based on the positive sum property of Graver
bases \cite{Hemmecke:PSP} needs to break the symmetry, see Section
\ref{Section: Computing the Graver Basis Faster}. Nonetheless, we
show how to exploit the symmetries also in this situation and
arrive at an even faster ``symmetric'' algorithm.

\begin{algorithm} \label{Slow Algorithm to Compute Graver Bases}
{(Algorithm to Compute $\GIP(\Lambda)$)}

\Kw{Input:} generating set $F$ of $\Lambda$ over $\Z$

\Kw{Output:} a set $G$ which contains $\GIP(\Lambda)$

\vspace{0.2cm} $G:=F\cup -F$

$C:=\bigcup\limits_{f,g\in G} \{f+g\}$

\Kw{while} $C\neq \emptyset $ \Kw{do}

\hspace{1.0cm}$s:=$ an element in $C$

\hspace{1.0cm}$C:=C\setminus\{s\}$

\hspace{1.0cm}$f:=\nF(s,G)$

\hspace{1.0cm}\Kw{if} $f\neq 0$ \Kw{then}

\hspace{2.0cm}$G:=G\cup \{f\}$

\hspace{2.0cm}$C:=C\cup\bigcup\limits_{g\in G} \{f+g\}$

\Kw{return} $G$.
\end{algorithm}

Behind the function $\nF(s,G)$ there is the following algorithm.
It aims at finding a relation $s=g_1+\ldots+g_r$, $g_i\in G$,
$g_i\red s$, $i=1,\ldots, r$. The function $\nF(s,G)$ returns $0$
if it found such vectors $g_i\in G$, or it returns a vector
$f\in\Lambda$ such that such a desired relation exists with
elements from $G\cup\{f\}$.

\begin{algorithm} \label{Slow NormalForm algorithm}
{(Normal Form Algorithm)}

\Kw{Input:} a vector $s$, a set $G$ of vectors

\Kw{Output:} a normal form of $s$ with respect to $G$

\vspace{0.2cm}

\Kw{while} there is some $g\in G$ such that $g\red s$ \Kw{do}

\hspace{1.0cm} $s:=s-g$

\Kw{return} $s$
\end{algorithm}

This algorithm due to Pottier always terminates and the set
$\GIP(\Lambda)$ is exactly the set of all $\red$-minimal vectors
in the final output $G$. Termination is guaranteed by the
following lemma, which we will employ again later.

\begin{lemma} \label{Gordan-Dickson Lemma, Sequence version}
{(Gordan-Dickson Lemma, Sequence version,
\cite{Cox+Little+OShea:92})}

Let $\{p_1,p_2,\ldots\}$ be a sequence of points in $\Z^n_+$ such
that $p_i\not\leq p_j$ whenever $i<j$. Then this sequence is
finite.
\end{lemma}

Now let us adapt the algorithm to exploit the symmetries. Let us
start with an immediate consequence of Lemma \ref{Orbits in G(A)
are always full}.

\begin{corollary} \label{Cor: Orbits in G(A) are always full}
If $s=g_1+\ldots+g_r$, $g_i\red s$, $i=1,\ldots,
r$, then $\sigma(s)=\sigma(g_1)+\ldots+\sigma(g_r)$,
$\sigma(g_i)\red\sigma(s)$, $i=1,\ldots, r$ for every
$\sigma\in\SA$. Thus, if $G=\bigcup\limits_{v\in
G^{\text{sym}}}\orb(v)$ and if $\nF(s,G)=0$ then
$\nF(\sigma(s),G)=0$ for every $\sigma\in\SA$.
\end{corollary}

In other words, if a representation $s=g_1+\ldots+g_r$, $g_i\in
G$, $g_i\red s$, $i=1,\ldots,r$ has been found, the symmetry of
$\Lambda$ already guarantees existence of a similar representation
for every element in $\orb(s)$.

Finally, we are in the position to exploit symmetries in Pottier's
algorithm. The main differences to the original algorithm will be
that instead of keeping the sets $G$ and $C$ in memory, we only
store their representatives under the given symmetry in sets
$G^{\text{sym}}$ and $C^{\text{sym}}$. (At any point during the
``symmetric'' algorithm we may go back to the original algorithm
by replacing all elements in $G^{\text{sym}}$ and $C^{\text{sym}}$
by the vectors from their orbits under $\SA$.)

Moreover, there a few more changes. Once we have found a nonzero
vector $f$ that is to be added to $G$ (and to $G^{\text{sym}}$ as
a new representative), we assume that we add the full orbit
$\orb(f)$ to $G$, as we know that the Graver basis would contain
the full orbit if $f$ was in fact a Graver basis element.
Accordingly, instead of adding only the vectors
\[
\bigcup\limits_{g\in G} \{f+g\}
\]
to $C$, we immediately include the vectors
\[
\bigcup\limits_{f'\in\orb(f),g\in G} \{f'+g\}.
\]
As $G$ will always be a union of full orbits, this last expression can
be transformed to
\[
\bigcup\limits_{f'\in\orb(f),g\in G} \{f'+g\}=
\bigcup\limits_{g\in G^{\text{sym}}}
\bigcup\limits_{\begin{array}{c}f'\in\orb(f)\\g'\in\orb(g)\end{array}}
\{f'+g'\}=\bigcup\limits_{g\in G^{\text{sym}}}
\bigcup\limits_{g'\in\orb(g)}\orb(f+g').
\]
Therefore, we update $C^{\text{sym}}$ as follows:
\[
C^{\text{sym}}=C^{\text{sym}}\cup\bigcup\limits_{g\in G^{\text{sym}}}
\bigcup\limits_{g'\in\orb(g)}\{f+g'\}.
\]
Note that since
\[
\bigcup\limits_{g'\in\orb(g)}\orb(f+g')=\bigcup\limits_{f'\in\orb(f)}\orb(f'+g)
\]
we have a choice in adding either all vectors $\{f+g'\}$ or all
vectors $\{f'+g\}$ to $C^{\text{sym}}$. Clearly, one would choose
to add as few new representatives to $C^{\text{sym}}$ as possible
to keep the number of S-vectors that need to be reduced small.
After all, these reductions are the most expensive part of the
algorithm. In the following algorithm, $\rep(H)$ for a set
$H\subseteq\Lambda$ of vectors shall denote a set of
representatives of $H$ under the symmetry group $\SA$.

\begin{algorithm} \label{Slow Algorithm to Compute Symmetric Graver Bases}
{(Algorithm to Compute $G^{\text{sym}}(\Lambda)$)}

\Kw{Input:} generating set $F$ of $\Lambda$ over $\Z$

\Kw{Output:} a set $G$ which contains $\GIP(\Lambda)$

\vspace{0.2cm} $G^{\text{sym}}:=\rep(F\cup -F)$
\hfill $G:=\orb(F\cup -F)$

\vspace{0.2cm} $C^{\text{sym}}:=\rep\left(\bigcup\limits_{f,g\in G}
\{f+g\}\right)$
\hfill $C:=\bigcup\limits_{f,g\in G} \{f+g\}$

\Kw{while} $C^{\text{sym}}\neq \emptyset $ \Kw{do}

\hspace{1.0cm}$s:=$ an element in $C^{\text{sym}}$

\hspace{1.0cm}$C^{\text{sym}}:=C^{\text{sym}}\setminus\{s\}$
\hfill $C:=C\setminus\orb(s)$

\hspace{1.0cm}$f:=\nF(s,G^{\text{sym}}):=\nF(s,G)$
\hfill $f:=\nF(s,G)$

\hspace{1.0cm}\Kw{if} $f\neq 0$ \Kw{then}

\hspace{2.0cm}$G^{\text{sym}}:=G^{\text{sym}}\cup\{f\}$ \hfill
$G:=G\cup\orb(f)$

\hspace{2.0cm}$C^{\text{sym}}:=C^{\text{sym}}\cup
\bigcup\limits_{g\in G} \bigcup\limits_{g'\in\orb(g)}\{f+g'\}$
\hfill $C:=C\cup\bigcup\limits_{f'\in\orb(f),g\in G} \{f'+g\}$

\Kw{return} $G=\orb(G^{\text{sym}})$.
\end{algorithm}

In this algorithm, we compute $\nF(s,G^{\text{sym}})$ via
$\nF(s,G^{\text{sym}}):=\nF(s,G)$. Clearly, from a practical
perspective, one would not want to keep the huge set $G$ in
memory. Then, of course, one needs to think about how to compute
$\nF(s,G^{\text{sym}})$ efficiently if only $G^{\text{sym}}$
instead of $G$ is available. (For example, $G$ might be simply too
big to be kept in memory.) This is still an open question and any
significant improvement in the solution of this problem would lead
to an equally significant improvement of the overall algorithm.

\begin{lemma}
Algorithm \ref{Slow Algorithm to Compute Symmetric Graver Bases}
always terminates and returns a set $G$ containing
$\GIP(\Lambda)$.
\end{lemma}

\boproof To prove termination, consider the sequence of vectors in
$G^{\text{sym}}\setminus\rep(F\cup -F)=\{f_1,f_2,\ldots\}$ as they are added to
$G^{\text{sym}}$ during the run of the algorithm. By construction, we
have $f_i\not\red f_j$, that is $(f_i^+,f_i^-)\not\leq (f_j^+,f_j^-)$
whenever $i<j$. Thus, by the Gordan-Dickson Lemma, this sequence must
be finite and the algorithm terminates.

Note that throughout the run of the algorithm, we always have
$G=\orb(G^{\text{sym}})$. Upon termination we know that
$\nF(f+g,G)=0$ for every pair of vectors $f,g\in G$. Thus, $G$
must contain the Graver basis $\GIP(\Lambda)$. \eoproof

\section{Computing the Graver basis faster}
\label{Section: Computing the Graver Basis Faster}

In this section we introduce a special generating set of
$\Lambda$. This set not only decreases the number of sums $f+g$
that need to be added to $C$, but more importantly, it reduces the
amount of work to compute $\nF(s,G)$ tremendously. As the latter
computation is the most expensive part of Pottier's algorithm,
this heavily speeds up the computation of $\GIP(\Lambda)$.
Moreover, we introduce a so-called critical-pair selection
strategy that chooses the next element $s$ from $C$ according to a
certain rule. This, together with the special input set, will
imply that the set $G$ returned by our algorithm is {\it exactly}
the Graver basis $\GIP(\Lambda)$. Not a single unnecessary vector
is computed!

In Pottier's algorithm, one has to wait until the very end to
extract the Graver basis from the returned set $G$. In contrast to
this, our algorithm provides a {\it certificate} of
$\red$-minimality for each vector that is added to $G$.
Consequently, at any point during the computation, a subset $G$ of
$\GIP(\Lambda)$ is known.

Let us assume from now on that $\Lambda$ is generated by $d$
vectors and that the first $d$ components of the vectors of
$\Lambda$ are linearly independent, that is, the only
$d$-dimensional vector orthogonal to the projection of $\Lambda$
to the first $d$ components is the zero vector. (This condition
can easily be achieved algorithmically on a lattice basis of
$\Lambda$ by integer row operations and by switching components.)
Let $\pi$ denote the projection of an $n$-dimensional vector onto
its first $d$ components. As the last $n-d$ components of
$\Lambda$ are linearly dependent on the first $d$ components, a
vector $\pi(v)\in\pi(\Lambda)$ can be uniquely lifted back to
$v\in\Lambda$. Moreover, this can easily be done algorithmically.

Next let us define a norm $\|.\|$ on vectors $v\in\Lambda$ as
follows: $\|v\|:=\|\pi(v)\|_1$, where $\|.\|_1$ denotes the
$L_1$-norm on $\R^d$. It can easily be checked that, under our
assumptions on $\Lambda$, this defines indeed a norm on $\Lambda$.
It is this norm definition that breaks existing symmetry of the
given problem.

Finally, we are ready to state the algorithm behind the implementation
in $\4ti2$ that seemingly defines the current state-of-the-art in the
computation of Graver bases.

Let $\pi(\bar{F})$ denote set of all $\red$-minimal nonzero
vectors in $\pi(\Lambda)$. Note that $\pi(\bar{F})$ has the
positive sum property with respect to $\pi(\Lambda)$, that is,
every vector $\pi(v)\in\pi(\Lambda)$ can be written as a positive
integer linear combination $\pi(v)=\sum\alpha_i\pi(f_i)$ of
vectors $\pi(f_i)$ from the set $\pi(\bar{F})$ that lie all in the
same orthant as $\pi(v)$. In other words, $\alpha_i\in\Z_{>0}$ and
$\pi(f_i)\red\pi(v)$ for all $i$ in this linear combination.
Clearly, this nice integer linear combination
$\pi(v)=\sum\alpha_i\pi(f_i)$ can be uniquely lifted to
$v=\sum\alpha_i f_i$ showing that $\bar{F}$ is in particular a
generating set of $\Lambda$ over $\Z$.

Moreover, this set $\pi(\bar{F})$ can be computed for example via
Pottier's algorithm once a lattice basis $F$ of
$\pi(\Lambda)\subseteq\Z^{d}$ over $\Z$ is given as input. Note
that all these $\red$-minimal vectors in $\pi(\bar{F})$ must lift
to $\red$-minimal elements in $\Lambda$, since they are already
indecomposable/minimal on the first $d$ components. Therefore, the
computation of $\pi(\bar{F})$ is usually far less expensive than
computing $\GIP(\Lambda)$ itself. In our computational
experiments, it often happened that this set $\pi(\bar{F})$ simply
contained the unit vectors in $\Z^{d}$ and their negatives.

\begin{algorithm} \label{Fast Algorithm to Compute Graver Bases}
{(Faster Algorithm to Compute $\GIP(\Lambda)$)}

\Kw{Input:} set $\bar{F}\subseteq\Lambda$ such that $\pi(\bar{F})$
is the set of $\red$-minimal nonzero vectors in $\pi(\Lambda)$

\Kw{Output:} a set $G$ which contains $\GIP(\Lambda)$

\vspace{0.2cm} $G:=\bar{F}$

$C:=\bigcup\limits_{f,g\in G: \pi(f) \text{ and } \pi(g) \text{ lie in
  the same orthant of } \R^{d}} \{f+g\}$

\Kw{while} $C\neq \emptyset $ \Kw{do}

\hspace{1.0cm}$s:=$ an element in $C$ with smallest $\|.\|$-norm

\hspace{1.0cm}$C:=C\setminus\{s\}$

\hspace{1.0cm}$f:=\nF(s,G)$

\hspace{1.0cm}\Kw{if} $f\neq 0$ \Kw{then}

\hspace{2.0cm}$G:=G\cup \{f\}$

\hspace{2.0cm}$C:=C\cup\bigcup\limits_{g\in G: \pi(f) \text{ and }
  \pi(g) \text{ lie in the same orthant of } \R^{d}} \{f+g\}$

\Kw{return} $G$.
\end{algorithm}

Due to our special input set, the function $\nF(s,G)$ can be sped
up as follows.

\newpage
\begin{algorithm} \label{Fast NormalForm algorithm}
{(Faster Normal Form Algorithm)}

\Kw{Input:} a vector $s$, a set $G$ of vectors

\Kw{Output:} a normal form of $s$ with respect to $G$

\vspace{0.2cm}

\Kw{if} there is some $g\in G$ such that $g\red s$ \Kw{return} $0$

\Kw{return} $s$
\end{algorithm}

Once a vector $g\in G$ with $g\red s$ is found, a representation
$s=g_1+\ldots+g_r$, $g_i\in G$, $g_i\red s$, $i=1,\ldots, r$, must
exist. Therefore, the function $\nF(s,G)$ can return $0$
immediately without explicitly constructing such a relation.

Since we started with a very special input set, only those pairs
of vectors $f,g\in G$ lead to a critical vector in $C$, for which
the projections $\pi(f)$ and $\pi(g)$ lie in the same orthant of
$\R^{d}$. Finally, let us prove our claims.

\begin{lemma}\label{Fast Algorithm to Compute Graver Bases is Correct}
Algorithm \ref{Fast Algorithm to Compute Graver Bases} always
terminates and returns a set $G$ containing $\GIP(\Lambda)$.
\end{lemma}

\boproof To prove termination, consider the sequence of vectors in
$G\setminus \bar{F}=\{f_1,f_2,\ldots\}$ as they are added to $G$
during the run of the algorithm. By construction, we have $f_i\not\red
f_j$, that is $(f_i^+,f_i^-)\not\leq (f_j^+,f_j^-)$ whenever
$i<j$. Thus, by the Gordan-Dickson Lemma, this sequence must be finite
and the algorithm terminates.

To prove correctness, let us assume that $z\in\GIP(\Lambda)$ is
not contained in the final output set $G$ of Algorithm \ref{Fast
Algorithm to Compute Graver Bases}. Without loss of generality we
may assume that $z$ has a smallest norm $\|z\|$ among all such
vectors from $\GIP(\Lambda)$. Therefore, we can assume that all
Graver basis elements $g$ with $\|g\|<\|z\|$ are contained in $G$.
In the following, we construct a contradiction to the assumption
$z\not\in G$ and correctness of Algorithm \ref{Fast Algorithm to
Compute Graver Bases} is proved.

As $\bar{F}$ is contained in $G$, there is a representation
$z=\sum\alpha_i v_i$ with positive integers $\alpha_i$ and vectors
$v_i\in G$ with $\pi(v_i)\red\pi(z)$. From the set of all such
linear integer combinations choose one such that
$\sum\alpha_i\|v_i\|_1$ is minimal.

Let us assume first that $\sum\alpha_i\|v_i\|_1>\|z\|_1$.
Therefore, there have to exist vectors $v_{i_1}, v_{i_2}$ in this
representation which have some component $k=k_0$ of different
signs. By construction, $k_0>d$, as the $v_i$ have all the same
sign as $z$ on the first $d$ components.

The vector $v_{i_1}+v_{i_2}$ was added to $C$ during the run of
Algorithm \ref{Fast Algorithm to Compute Graver Bases} (as
$\pi(v_{i_1})$ and $\pi(v_{i_2})$ lie in the same orthant of
$\R^{d}$ by construction). If $\|v_{i_1}+v_{i_2}\|=\|z\|$, then
all other $v_i$, $i\neq i_1,i_2$, must satisfy $\|v_i\|=0$ as
$\pi(v_i)\red\pi(z)$ for all $i$. But $\pi(v_i)=0$ implies $v_i=0$
and thus $v_{i_1}+v_{i_2}=z$. Since $v_{i_1}+v_{i_2}(=z)$ is a
vector that was added to $C$ during the run of the algorithm, the
vector $z$ is eventually chosen as $s\in C$. Being $\red$-minimal,
we have $\nF(z,G)=z$ and thus, $z$ must have been added to $G$, in
contradiction to our assumption $z\not\in G$.

Therefore, we may assume that $\|v_{i_1}+v_{i_2}\|<\|z\|$.
However, since all Graver basis elements $v$ with norm
$\|v\|<\|z\|$ are assumed to be in $G$, there must exist a
representation $v_{i_1}+v_{i_2}=\sum\beta_j v'_j$ for finitely
many $\beta_j\in\Z_{>0}$, $v'_j\in G$, and $\beta_j v'_j\red
v_{i_1}+v_{i_2}$ for all $j$. This implies that we have for each
component $k=1,\ldots,n$,
\[
\sum_j\beta_j |{v'_j}^{(k)}|=|\sum_j\beta_j {v'_j}^{(k)}|
= |(v_{i_1}+v_{i_2})^{(k)}|\leq |v_{i_1}^{(k)}|+|v_{i_2}^{(k)}|,
\]
where the last inequality is strict for $k=k_0$ by construction.
Summing up over $k=1,\ldots,n$, yields $\sum \beta_j
\|{v'_j}\|_1=\|v_{i_1}+v_{i_2}\|_1<\|v_{i_1}\|_1+\|v_{i_2}\|_1$.
But now $z$ can be represented as
\begin{eqnarray*}
z & = & \alpha_{i_1}v_{i_1}+\alpha_{i_2}v_{i_2}
+\sum_{i\neq i_1,i_2}{\alpha_i v_i}\\
  & = & \sum\beta_j v'_j+(\alpha_{i_1}-1)v_{i_1}+(\alpha_{i_2}-1)v_{i_2}
+\sum_{i\neq i_1,i_2}{\alpha_i v_i}
\end{eqnarray*}
and it holds
\[
\sum\beta_j \|v'_j\|_1+(\alpha_{i_1}-1)\|v_{i_1}\|_1+
(\alpha_{i_2}-1)\|v_{i_2}\|_1+\sum_{i\neq i_1,i_2}{\alpha_i\|v_i\|_1}<
\sum\alpha_i\|v_i\|_1
\]
in contradiction to the minimality required on
$\sum\alpha_i\|v_i\|_1$. Thus, our assumption
$\sum\alpha_i\|v_i\|_1>\|z\|_1$ was wrong and
$\sum\alpha_i\|v_i\|_1=\|z\|_1$ must hold.

But this last equation implies that $v_i\red z$ for all $i$,
contradicting $\red$-minimality of $z$ unless the representation
$z=\sum\alpha_i v_i$ is trivial, that is $z=v_1\in G$. This,
however, again contradicts our initial assumption $z\not\in G$
and thus $\GIP(\Lambda)\subseteq G$. \eoproof

It should be noted that we did not make use of our selection strategy
to prove termination and correctness of Algorithm
\ref{Fast Algorithm to Compute Graver Bases}. It is the following
Lemma that provides a certificate for $\red$-minimality of vectors in
$G$.

\begin{lemma}\label{G=GIP}
The set $G$ returned by Algorithm \ref{Fast Algorithm to Compute
Graver Bases} equals $\GIP(\Lambda)$.
\end{lemma}

\boproof The main observation needed in this proof is that the
norms $\|s\|$ of the vectors $s$ that are chosen from $C$ form a
non-decreasing sequence. This follows from the definition that
$\nF(s,g)$ only returns either $0$ or $s$ and from our condition
that only vectors $f+g$ are added to $C$ whose first $d$
components have the same sign pattern. The latter implies
$\|f+g\|=\|f\|+\|g\|$.

Now assume that some $z\in G$ is not contained in $\GIP(\Lambda)$.
Thus, there is some $g\in\GIP(\Lambda)$ with $g\red z$. Under our
assumptions on $\Lambda$, $g\red z$ implies $\|g\|<\|s\|$. Since
$G$ contains $\GIP(\Lambda)$ (and thus in particular the vector
$g$) at the end of the algorithm and since after $z$ only vectors
$f$ are added to $G$ that have a norm $\|f\|\geq\|z\|$, the vector
$g$ must have been contained in $G$ already at the time when $z$
was added to $G$. This however, implies that Algorithm \ref{Fast
Algorithm to Compute Graver Bases} must have computed
$\nF(z,G)=0$, a contradiction to the assumption that $z$ was added
to $G$. \eoproof

\section{Computing the symmetric Graver basis faster}
\label{Section: Computing the symmetric Graver basis faster}

Although the definition of $\|.\|$ in the previous section broke
most if not all existing symmetry in the problem, we will now
combine the ideas of Sections \ref{Section: Computing the Graver
Basis} and \ref{Section: Computing the Graver Basis Faster} to a
faster algorithm to compute symmetric Graver bases. The main idea
is to use the norm $\|.\|$ defined on $\Z^n$ to define a norm on
orbits: For $T\subseteq\Lambda$, we define
$\|T\|:=\min\{\|v\|:v\in T\}$. Then the new symmetric Graver basis
algorithm looks as follows.

\begin{algorithm} \label{Fast Algorithm to Compute Symmetric Graver Bases}
{(Faster Algorithm to Compute $G^{\text{sym}}(\Lambda)$)}

\Kw{Input:} set $\bar{F}\subseteq\Lambda$ such that $\pi(\bar{F})$
is the set of $\red$-minimal nonzero vectors in $\pi(\Lambda)$

\Kw{Output:} a set $G$ which contains $\GIP(\Lambda)$

\vspace{0.2cm} $G^{\text{sym}}:=\rep(\bar{F})$

\vspace{0.2cm} $C^{\text{sym}}:=\rep\left(\bigcup\limits_{f,g\in G}
\{f+g\}\right)$

\Kw{while} $C^{\text{sym}}\neq \emptyset $ \Kw{do}

\hspace{1.0cm}$s:=$ an element in $C^{\text{sym}}$ with smallest value
of $\|\orb(s)\|$

\hspace{1.0cm}$C^{\text{sym}}:=C^{\text{sym}}\setminus\{s\}$

\hspace{1.0cm}$f:=\nF(s,G^{\text{sym}})$

\hspace{1.0cm}\Kw{if} $f\neq 0$ \Kw{then}

\hspace{2.0cm}$G^{\text{sym}}:=G^{\text{sym}}\cup\{f\}$

\hspace{2.0cm}$C^{\text{sym}}:=C^{\text{sym}}\cup
\bigcup\limits_{g\in G} \bigcup\limits_{g'\in\orb(g)}\{f+g'\}$

\Kw{return} $G=\orb(G^{\text{sym}})$.
\end{algorithm}

Due to our special input set and our norm defined on orbits, we
can again simplify and speed up the normal form computation
$\nF(s,G^{\text{sym}}):=\nF(s,\orb(G^{\text{sym}}))$ by using
Algorithm \ref{Fast NormalForm algorithm} instead of Algorithm
\ref{Slow NormalForm algorithm}. Again, from a practical
perspective, one would want to compute $\nF(s,G^{\text{sym}})$
without recovering or storing the huge set $\orb(G^{\text{sym}})$.

On the other hand, it should be noted that we do not, as in
Algorithm \ref{Fast Algorithm to Compute Graver Bases}, have an
orthant condition in Algorithm \ref{Fast Algorithm to Compute
Symmetric Graver Bases} that reduces the number of vectors added
to $C^{\text{sym}}$. This is done to ``undo'' the symmetry
breaking caused by $\|.\|$.

\begin{lemma}
\label{Fast Algorithm to Compute Symmetric Graver Bases is
Correct} Algorithm \ref{Fast Algorithm to Compute Symmetric Graver
Bases} always terminates and returns the set $G=\GIP(\Lambda)$.
\end{lemma}

\boproof The subsequent proof follows similar lines as the proof
of Lemma \ref{Fast Algorithm to Compute Graver Bases is Correct}.

To prove termination, we consider again the sequence of vectors in
$G^{\text{sym}}\setminus\rep(\bar{F})=\{f_1,f_2,\ldots\}$ as they
are added to $G^{\text{sym}}$ during the run of the algorithm. By
construction, we have $f_i\not\red f_j$, that is
$(f_i^+,f_i^-)\not\leq (f_j^+,f_j^-)$ whenever $i<j$. Thus, by the
Gordan-Dickson Lemma, this sequence must be finite and the
algorithm terminates.

Next we show that $\GIP(\Lambda)\subseteq G$. Assume on the
contrary that this is not the case. Among all elements
$z\in\GIP(\Lambda)\setminus G$ choose one with $\|\orb(z)\|$
smallest. Moreover, we may assume that $z$ is a representative of
$\orb(z)$ with $\|\pi(z)\|_1=\|\orb(z)\|$. By our generating
assumption on the set $\bar{F}\subseteq G$, there exists a
non-trivial representation $z=\sum\alpha_i v_i$ with positive
integers $\alpha_i$ and vectors $v_i\in G$ with
$\pi(v_i)\red\pi(z)$. From the set of all such linear integer
combinations choose one such that $\sum\alpha_i\|v_i\|_1$ is
minimal. Note that from $\pi(v_i)\red\pi(z)$ and the fact that the
relation $z=\sum\alpha_i v_i$ is non-trivial, we conclude
$\|\pi(v_i)\|_1<\|\pi(z)\|_1$ and thus $\|\orb(v_i)\|<\|\orb(z)\|$
for all $i$.

Let us assume first that $\sum\alpha_i\|v_i\|_1>\|z\|_1$.
Therefore, there have to exist vectors $v_{i_1}, v_{i_2}$ in this
representation which have some component $k=k_0$ of different
signs. By construction, $k_0>d$, as the $v_i$ have all the same
sign as $z$ on the first $d$ components.

The orbit of the vector $v_{i_1}+v_{i_2}$ was added to
$C^{\text{sym}}$ during the run of the algorithm. If
$\|v_{i_1}+v_{i_2}\|=\|z\|$, then all other $v_i$, $i\neq i_1,i_2$
must satisfy $\|v_i\|=0$ as $\pi(v_i)\red\pi(z)$ for all $i$. But
$\pi(v_i)=0$ implies $v_i=0$ and thus $v_{i_1}+v_{i_2}=z$. Since
$v_{i_1}+v_{i_2}(=z)$ is a vector whose orbit (or better: a
representative of it) was added to $C^{\text{sym}}$ during the run
of the algorithm, a representative of the orbit of the vector $z$
is eventually chosen as $s\in C^{\text{sym}}$. Being
$\red$-minimal, we have
$\nF(\rep(\orb(z)),G^{\text{sym}})=\rep(\orb(z))$ and thus,
$\rep(\orb(z))$ must have been added to $G^{\text{sym}}$, in
contradiction to our assumption $z\not\in G$.

Therefore, we may assume that $\|v_{i_1}+v_{i_2}\|<\|z\|$.
However, since all Graver basis elements with norm strictly
smaller than $\|z\|$ are assumed to be in $G$, we may continue
literally as in the proof of Lemma \ref{Fast Algorithm to Compute
Graver Bases is Correct}: rewrite $z=\sum\alpha_i v_i$ and arrive
at a contradiction to the minimality of $\sum\alpha_i v_i$.

Thus, we must have $\sum\alpha_i\|v_i\|_1=\|z\|_1$, which implies
that $v_i\red z$ for all $i$. As $z$ is $\red$-minimal, this is
only possible if the relation $z=\sum\alpha_i v_i$ is trivial,
that is $z=v_1$. As $v_1\in G$ we now also have $z\in G$ as
desired.

To show that not only $\GIP(\Lambda)\subseteq G$ but in fact
$\GIP(\Lambda)=G$ is true, assume
$G\setminus\GIP(\Lambda)\neq\emptyset$. Among all elements $z\in
G\setminus\GIP(\Lambda)$ choose one with $\|\orb(z)\|$ smallest.
Moreover, we may again assume that $z$ is a representative of
$\orb(z)$ with $\|\pi(z)\|_1=\|\orb(z)\|$. Note that by Lemma
\ref{Orbits in G(A) are always full} and since
$z\not\in\GIP(\Lambda)$, not a single element in $\orb(z)$ belongs
to $\GIP(\Lambda)$.

We now start with a non-trivial representation $z=\sum\alpha_i
v_i$ with positive integers $\alpha_i$ and vectors
$v_i\in\bar{F}\subseteq G$ with $\pi(v_i)\red\pi(z)$. Clearly,
$\pi(v_i)\red\pi(z)$ implies $\|\pi(v_i)\|_1<\|\pi(z)\|_1$ and
thus $\|\orb(v_i)\|<\|\orb(z)\|$ for all $i$. Next, we follow
similar steps as above (or more precisely: as in the proof of
Lemma \ref{Fast Algorithm to Compute Graver Bases is Correct}) to
change $z=\sum\alpha_i v_i$ into a representation $z=\sum\beta_j
w_j$ with $\beta_j>0$, $w_j\in G$ and $w_j\red z$ for all $j$. As
$z$ is not $\red$-minimal, this representation is not trivial.

Note that in $z=\sum\alpha_i v_i$, $z=\sum\beta_j w_j$, and in any
intermediate representation $z=\sum\beta_k u_k$, we have that the
summands $v_i,w_j,u_k$ have a strictly smaller norm on the first
$d$ components as $z$. In fact, the same property holds for the
sum of two summands that are iteratively needed for the rewriting
steps. Thus, their corresponding orbits also have a strictly
smaller norm than $\|\orb(z)\|$. We conclude that at the time
$\rep(\orb(z))$ was chosen from $C^{\text{sym}}$ and then added to
$G^{\text{sym}}$, representatives for all orbits $\orb(w_j)$ were
already in $G^{\text{sym}}$. Thus,
$\nF(\rep(\orb(z)),G^{\text{sym}})$ must have returned $0$ in
contradiction to the assumption that $\rep(\orb(z))$ was added to
$G^{\text{sym}}$.

We conclude $G\setminus\GIP(\Lambda)=\emptyset$ and hence
$G=\GIP(\Lambda)$. \eoproof

\section{Computational Experiments}

In this section we report on computational experience with a few
examples that have symmetry. These problem deal all with $3$-way
tables, that is, with $k_1\times k_2\times k_3$ tables of unit
cubes, each containing an integer number. For each problem, the
lattice $\Lambda$ is the set of lattice points in the kernel of
the matrix that encodes the conditions that the sums along each
one-dimensional row parallel to a coordinate axis are $0$. Example
\ref{Example: 3x3 tables} presented the case of $3\times 3$
tables.

The following running times demonstrate that, as expected, our
symmetric algorithm heavily speeds up the computation. (Times are
given in seconds on a Sun UltraSparc III+ with 1.05 GHz.)

\begin{center}
\begin{tabular}{|l|r|r|r|r|r|}
\hline Problem & $|\SA|$ & Size of Graver basis &
$|G^{\text{sym}}|$ & Algorithm~\ref{Fast Algorithm to Compute
Graver Bases} &
Algorithm~\ref{Fast Algorithm to Compute Symmetric Graver Bases} \\
\hline
$3\times 3\times 3$ & $1,296$ &       $795$ &   $7$ &       $2$ &      $1$ \\
$3\times 3\times 4$ & $1,728$ &    $19,722$ &  $27$ &   $1,176$ &      $9$ \\
$3\times 3\times 5$ & $8,640$ &   $263,610$ &  $61$ & $560,517$ &    $526$ \\
$3\times 4\times 4$ & $6,912$ & $4,617,444$ & $784$ &      $--$& $260,590$ \\
\hline
\end{tabular}
\end{center}

\end{document}